\documentclass{article}

\usepackage{PRIMEarxiv}

\usepackage[utf8]{inputenc} 
\usepackage[T1]{fontenc}    
\usepackage{hyperref}       
\usepackage{url}            
\usepackage{booktabs}       
\usepackage{amsfonts}       
\usepackage{nicefrac}       
\usepackage{microtype}      
\usepackage{lipsum}
\usepackage{graphicx}
\graphicspath{{media/}}     
\usepackage{amsmath}
  \usepackage{amssymb}
  \newtheorem{theorem}{Theorem}[section]

\newtheorem{remark}[theorem]{Remark}
\DeclareUnicodeCharacter{2212}{-}
\usepackage{verbatim}
\usepackage{xspace}
\usepackage{amssymb,enumitem}
\title{Smoothness and analyticity of $f'=\exp({f^{-1}})$
}

\author{
  Zeraoulia Rafik \\
  Affiliation \\
  University of batna2.Algeria \\
  City Batna,Hamla3\\
    \texttt{\{ Email\}r.zeraoulia@univ-batna2.dz} \\
}

\begin{document}
\maketitle

\begin{abstract}
This paper considers some analytical and numerical aspects of the problem defined by an equation  of the type $f'=\exp({f^{-1}})$ with $f^{-1}$ is a composional inverse of $f$,some new  analytical and numerical results are presented using RK4 and Explicit Runge kuta methods. 
\end{abstract}


\section{Introduction}
During the last three decades, the participants of the Perm Seminar \cite{AZB:03} developed a unified theory of a broad generalization of the differential equation. The conducted studies have established a close relationship between numerous tasks previously studied outside
connections with each other, and made it possible to propose more advanced methods for solving them.

The basic elements of the modern theory of functional differential equations were established in \cite{A.D:49}, \cite{Myshkis:55}. The years after \cite{Myshkis:55} have witnessed an explosive development of the theory of functional differential equations and their applications (see, e.g. \cite{Maksimov:66} , \cite{Rak:2004}, \cite{J:71}, \cite{V.B:86}, \cite{V.kol:92})and the numerous references therein).

A functional differential equation (also called a differential equation with deviating argument, cf. also Differential equations, ordinary, with distributed arguments) can be considered as a combination of differential and functional equations. The values of the argument in a functional differential equation can be discrete, continuous or mixed. Correspondingly, one introduce the notions of a differential difference equation, an integro differential equation, etc.
Let $\mathbb{K}$. be an algebraically closed field of characteristic zero. An algebraic ordinary
differential equation (AODE) is an equation of the form:
\begin{equation*}
    F\bigg(x,y,\frac{dy}{dx} \cdots \frac{d^n y}{dx^n}\bigg)=0
\end{equation*}
for some $n\in \mathbb{N}$ and $F$ a polynomial in $n+2$ variables over $\mathbb{K}$. This paper addresses analyticity and functions satisfying  nonlinear AODS 
\begin{equation}\label{1}
f'=\exp(f^{-1})
\end{equation}
with $f^{-1}$ which map $\mathbb{R}$ to $\mathbb{R}$,the solution of (\ref{1}) is a formel power series arround $0$ satisfying $f(0)=0$ ,In mathematics, formal power series are a generalization of polynomials as formal objects, where the number of
terms is allowed to be infinite, this implies giving up the possibility to substitute arbitrary values for
indeterminates. This perspective contrasts with that of power series, whose variables designate numerical values,
and which series therefore only have a definite value if convergence can be established. Formal power series are
often used merely to represent the whole collection of their coefficients. In combinatorics, they provide
representations of numerical sequences and of multisets, and for instance allow giving concise expressions for
recursively defined sequences regardless of whether the recursion can be explicitly solved , this is known as the
method of generating functions.
The problem of finding formal power series solutions of AODEs has a long
history and it has been heavily studied in the literature. The Newton polygon
method is a well known method developed for studying this problem. In \cite{C.Briot:56},
Briot and Bouquet use the Newton polygon method for studying the singularities
of first order and first degree ODEs. Fine gave a generalization of the method
for arbitrary order AODEs in \cite{H.B:89}. By using Newton polygon method, one
can obtain interesting results on a larger class of series solutions which is called
generalized formal power series solutions, i.e. power series with real exponents. In
\cite{D:91} Grigoriev and Singer proposed a parametric version of the Newton polygon
method and use it to study generalized formal power series solutions of AODEs.
A worth interpretation of the parametric Newton polygon method can be found in
\cite{P.Forty:09} and \cite{Cano:05}. However, it has been shown in \cite{Dora:97} that Newton polygon
method for AODEs has its own limits in the sense that in some cases, it fails to give a solution.Formal power series allow one to employ much of the analytical machinery of power series in settings which don't have natural notions of convergence. They are also useful in order to compactly describe sequences and to find closed formulas for recursively described sequences , this is known as the method of generating functions ,One application of generating functions is to solve counting problems like with Rook Polynomials  \cite{Feryal:13} . In this setting the coefficients will count a slight variant of the general problem and there are often infinitely many. This is useful because we can often calculate the coefficients in isolation without having to calculate the other terms to solve these problems.As we claimed before the aim of this paper is looking analyticity and behavior of the formel power series the solution of (\ref{1}) some new results regarding the functional differential equation defined in (\ref{1}) are presented.

\section{Main results}
\begin{enumerate}
   \item   There is no such function $f:\mathbb{R}\to \mathbb{R}$ satisfying
$f'=\exp(f^{-1})$ (there is no analytical solution at $0$ to $f'=\exp(f^{-1})$ with $f(0)=0$)

\item  There is a unique solution in formal power series arround $0$  satisfying $f(0)=0$

\item  The general term of the obtained formel power series $b_n=(-1)^{n}n!a_{n}$ appears  to always be a positive integer for $n >1$ (partial result)

\item  Both $f_{-}$ and $f_{+}$ are smooth functions, and they form an orbit of order at most $2$ of the Picard iteration of  $f_{n+1}'= \exp({f_n^{-1}})$ with the initial conditions $f_{n+1}(0)=0$ and they have the same Taylor expansion as calculated using formal power series expansion.

\end{enumerate}

For the proof of the first result  we have noticed that there is no such function. Since $f$ would have to map $\mathbb{R}$  onto $\mathbb{R}$ for the equation (\ref{1}) to make sense at all $x\in\mathbb{R}$, it follows that $f^{-1}(x)\to -\infty$ also as $x\to -\infty$, so $f'\to 0$. Thus $f(x)\ge x$, say, for all small enough $x$, hence $f^{-1}(x)\le x$ eventually, but then the equation shows that $f'\le e^x$, which is integrable on $(-\infty, 0)$, so $f$ would approach a limit as $x\to -\infty$ and not be surjective after all and we are done.Now for the second result we may need to differentiate equation (\ref{1}) another time to get 
\begin{equation}\label{eq2}
   f''(x)=f'(x)\cdot (f^{-1})'(x). 
\end{equation}
or equivalently 
\begin{equation}\label{eq3}
    f''(x)\cdot f'(f^{-1}(x))=f'(x)
\end{equation}
We note that  simplification of equation (\ref{1}) to  (\ref{eq2}) or  (\ref{eq3}) looks less scary to us than having  $f^{-}$ as an exponent .We have calculated first few coefficients  of the unique power series solution using mathematica and we have got :
\begin{align*}
[0, 1, 1/2, 0, 1/24,&\\ -1/20, 13/180, -197/1680 &\\, 2101/10080,&\\ -48203/120960, 2938057/3628800, -23059441/13305600,&\\ 74408941/19160064,&\\ -9409883317/1037836800]\\
\end{align*}
here is the mathematica code ,one can refer to it for more coefficients .(just changing n in the code)

\begin{verbatim}
    n = 12;
Input
ansatz /.
First[SolveAlways[D[ansatz, {x, 2}]
 ComposeSeries[D[ansatz, x], InverseSeries[ansatz]] - D[ansatz, x] == 0, x]] 
                           
output
x + x^2/2 + x^4/24 - x^5/20 + (13 x^6)/180 - (197 x^7)/1680 +
   (2101 x^8)/10080 - (48203 x^9)/120960 + (2938057 x^10)/3628800 -
   (23059441 x^11)/13305600 + (74408941 x^12)/19160064 + O[x]^13
\end{verbatim}

The Calcualtion of $100$ first terms (coefficients) using the above mathematica code show us that is pretty clear that $|a_n|^{-1/n}\to 0$ as $n\to\infty$ , so the radius of convergence is zero, so this approach will not give a solution that is an actual function. It is also of some interest that the number  $b_n=(-1)^nn!a_n$ appears to always be a positive integer (for $n>1$), but this sequence is not in OEIS. Also $b_n$ it does not factorise in a way that suggests that there could be a simple formula: for example $b_{10}=2938057$, which is prime, just a simple modification in the above code,taking $C(k)\to (-1)^{k} b_{k} k! C(k)$ then run and see what happen .A formal Taylor series (e.g.f.) solution about the origin can be obtained a few ways.

Let $f^{(-1)}(x) = e^{b.x}$ with $(b.)^n=b_n \;$ and $ \; b_0=0$.
Then A036040  (Bell polynomials) gives the e.g.f (exponential generating function) \cite{F:18}. 
$$e^{f^{(-1)}(x)}= e^{e^{b.x}}= 1 + b_1 x + (b_2+b_1^2) \frac{x^2}{2!}+(b_3+3b_1b_2+b_1^3)\frac{x^3}{3!}+\cdots \; ,$$
and the Lagrange inversion / series reversion formula (LIF) A134685 gives
$$f'(x)= \frac{1}{b_1} + \frac{1}{b_1^3} (-b_2) x + \frac{1}{b_1^5}(3b_2^2-b_1b_3)\frac{x^2}{2!}+\cdots \; .$$
Equating the two series and solving recursively gives
$$b_n \rightarrow (0,1,-1,3,-16,126,-1333,...)$$
which is signed A214645. This follows from the application of the inverse function theorem (essentially the LIF again) ,one can refer to the \textbf{Theorem 3.6},page 15 in  \cite{Alan:09} and \cite{M:10}

$$f'(z) = 1/f^{(-1)}{'}(\omega) \; ,$$
when $(z,\omega)=(f^{(-1)}(\omega),f(z)) $, leading to
$$f^{(-1)}{'}(x) = \exp[-f^{(-1)}(f^{(-1)}(x))],$$
the differential equation defining signed A214645.
Applying the LIF to the sequence for $b_n$ gives the e.g.f. $f(x)=e^{a.x}$ equivalent of F.C.'s o.g.f.(Ordinary generating function) \cite{F:18}
$$ a_n \rightarrow (0,1,1,0,1,-6,52,...).$$
As another consistency check, apply the formalism of A133314 for finding the multiplicative inverse of an e.g.f. to find the e.g.f. for $\exp[-A(-x)]=\exp[f^{(-1)}(x)]$ from that for 
$$\exp[A(-x)]= 1 - x + 2 \frac{x^2}{2!}-7 \frac{x^3}{3!}+\cdots \; ,$$
which is signed A233335, as noted in A214645. This gives $f'(x)=a. \; e^{a.x}$. The inverse function theorem here might be more aptly called the inverse formal series theorem. As we can see, the differential equations and inverses here in analytic guise are concise statements of relations among the coefficients of formal series (e.g.f.s or o.g.f.s).(exponential generating formel series).
Now for the proof of last main  result ,The proof of the first result  exploits the growth of $f(x)$ when $x\to -\infty$ and we obtained a contradiction, which resolves the question of analyticity arround $0$ nicely with $f(0)=0$, but also invites the following question: what if we restrict to $f:\mathbb R_{\ge0}\to\mathbb R_{\ge0}$ and impose $f(0)=0$? This idea has been explored in the first proof, where a formal power series expansion is obtained for $f$ which does not seem to converge for any $x\ne0$.

Taking another approach, we can use an iteration scheme starting from $f_1(x)=x$ and inductively solve the ODE $f_{n+1}'=e^{f_n^{-1}}$ with the initial condition $f_{n+1}(0)=0$ to obtain $f_{n+1}$, much in the spirit of Picard iteration. Explicitly, for example, we have
$f_2'=e^x$ and $f_2=e^x-1$;
$f_3'=e^{\ln(x+1)}=1+x$ and $f_3=x+x^2/2$;
$f_4'=e^{\sqrt{1+2x}-1}$ and $f_4=e^{\sqrt{1+2x}-1}(\sqrt{1+2x}-1)$
and the next iteration produces non elementary functions. It is clear that the sequence $(f_{2k-1})_{k\ge1}$ is increasing, $(f_{2k})_{k\ge1}$ is decreasing, and $f_{2k-1}<f_{2k}$, so there are respective limits $f_-=\lim_{k\to\infty} f_{2k-1}$ and $f_+=\lim_{k\to\infty} f_{2k+1}$, with $f_-\le f_+$. It is also clear that from $n\ge2$ on the function $f_n'=e^{f_{n-1}^{-1}}$ is positive and increasing, so $f_n$ is increasing and convex, which can be passed to the limit to show that both $f_-$ and $f_+$ are also increasing and convex. As such they are continuous, and by Dinis theorem \cite{WU:20} $f_{2k-1}$ converges to $f_-$ locally uniformly and similarly for $f_+$. Furthermore, the inequality $|x-y|\le |f_n(x)-f_n(y)|$ (as $f_n'=e^{f_{n-1}^{-1}}\ge1$) can also be passed to the limit. Then the following chain of inequalities:
$$|f_-^{-1}(x)-f_{2k-1}^{-1}(x)|\le |x-f_-(f_{2k-1}^{-1}(x))|=|f_{2k-1}(f_{2k-1}^{-1}(x))-f_-(f_{2k-1}^{-1}(x))|$$
shows that $f_{2k-1}^{-1}$ converges locally uniformly to $f_-^{-1}$, which then implies $f_{2k}'$ converges locally uniformly to $e^{f_-^{-1}}$. Hence $f_+'=e^{f_-^{-1}}$, and similarly $f_-'=e^{f_-^{-1}}$. From this it can be shown that $f_{2k-1}$ converges to $f_-$ locally in $C^\infty$, so both $f_-$ and $f_+$ are smooth functions, and they form an orbit of order at most $2$ of the above iteration scheme. Moreover it can be shown that the first $n$ terms of the Taylor expansion of $f_n$ agrees with what have been calculated formally using the above mathematica code, so both $f_-$ and $f_+$ have the same Taylor expansion as calculated using formal power series expansion.

In light of the above, a priori the following three scenarios can happen:
\begin{itemize}
 \item1) $f_-\neq f_+$ and we have a genuine orbit of order 2, consisting of two functions having the same Taylor expansion at $0$ but not being identical.
 \item 2) $f_-=f_+$ is an actual solution to the equation $f'=e^{f^{-1}}$, but it is merely $C^\infty$ but not analytic, having a divergent power series expansion at $0$.
 \item3) $f_-=f_+$ is an actual solution to the equation $f'=e^{f^{-1}}$, and it is analytic on a neighborhood of $0$; we are just misled by the first $100$ or so terms of the Taylor expansion.

Now finally comes the question: which of the above scenario is the reality? In the first two scenarios, one can also ask what is the growth rate of $f_-(x)$ and $f_+(x)$ as $x\to +\infty$.
\end{itemize}
For the good and correct scenario choice of the obove scenarios we note that the convergence is not hard to demonstrate. For instance,if $f,g$ are two increasing functions with $f(0)=g(0)=0$ and $f(x),g(x)\geq x$ for all $x\geq 0$, and $F,G$ are their images under the Picard map, then for every $T>0$, the functional $\Phi(f,g,T)=\int_0^T|f(t)-g(t)|\,dt$ satisfies
\begin{equation}\label{picard}
\Phi(F,G,T)\le \int_0^T e^t\Phi(f,g,t)\,dt
\end{equation}
and it follows that on every finite interval $[0,T]$  we have convergence in  $L^1$ and therefore in $C^{\infty}$.

There is no analytic local solution at $0$ to $f'=e^{f^{-1}}$, $f(0)=0$, that is, the formal power series solution is diverging. Together with the solution given using the Picard map where the inequality (\ref{picard}) statisfied, this means the actual scenario is \textbf{2} For convenience of notation, We shall consider the equivalent equation 
\begin{equation}\label{eq an}
\begin{cases} 
g' =e^{g\circ g},  \\ g(0)=0, 
 \end{cases}
 \end{equation}
satisfied by $g(x):=-f^{-1}(-x)$ (Indeed, by the rule of the derivative of an inverse, $(f^{-1})'(x)={1\over f'(f^{-1}(x))}=e^{-f^{-1}(f^{-1}(x))}$ so that $g'(x) =e^{g(g(x))}$; (see the proof above of the second result)
Indeed, assume by contradiction the formal power series solution $x+{1\over2}x^2+{1\over2}x^3+{2\over3}x^4+\&c.$ to the above equation (\ref{eq an}) has a positive radius of convergence. Then, it extends uniquely by analytic continuation to a maximally defined analytic function, still denoted $g$ (that is, defined on the largest positive interval $[0,a)$, for some $0<a\le+\infty$). 

Note that the Taylor series of $g$ at $0$ has non negative coefficients. This follows immediately by induction, equating the  coefficients of $g'$ and $e^{g\circ g}$; incidentally, this series is the EGF of the positive integer sequence  A214645,  As a consequence (one can check the details below), $g$ is totally monotonic on $[0,a)$; in particular $g'(x)>g'(0)=1$ and $g(x)>x$ for all $0<x<a$, and $g$ is invertible. 

Then we observe that $\log( g'( g^{-1}(x))$  is a well defined analytic function on the interval $g[0,a)$, and coincides with $g$ locally at $0$. By the maximality of $[0,a)$ we have thus $g[0,a)\subset[0,a)$, but, due to the inequality $g(x)>x$ on $(0,a)$, this inclusion is only possible if $a=+\infty$, so that $g$ is unbounded. On the other hand,since $e^{-g(g(t))}g'(t)=1$ and $g(t)\ge t$, we have for any $x\ge0$
$$x=\int_0^{x}e^{-g(g(t))}g'(t)dt=\int_0^{g(x)}e^{-g(s)}ds\le \int_0^{+\infty}e^{-s}ds=1 ,$$
a contradiction. 
\begin{remark}
 To justify the total monotonicity of $g$, note that, as a general elementary fact, a  real analytic function on an interval $I$, whose Taylor series at some point $x_0\in I$ has non negative coefficients, has Taylor series with non negative coefficients ay any point $x\in I$, $x\ge x_0$.* Indeed, this is clear for $x_1\ge x_0$ within the radius of convergence of $x_0$, and since there is a uniform radius of convergence at any $y\in [x_0,x]$, one reaches $x$ by finitely many steps $x_0<x_1<\dots<x_n=x$. 
 In fact more is true: a real analytic function on $\mathbb R$, whose Taylor series at some point $x_0\in\mathbb R$ has non negative coefficients is an entire function, so that any $x\ge x_0$ is reached in just one step. 
\end{remark}
\begin{remark}
The very same argument works for other differential functional equations like e.g.
$$\begin{cases} g' =1 + {g\circ g},  \\ g(0)=0, 
 \end{cases}$$
that generates the sequence OEIS A001028. 
As before, a maximally defined analytic solution $g$, if any, must be totally monotonic and defined for all  $x\ge0$, for otherwise $ g'\circ g^{-1} -1$ would be a proper extension of it. Then we reach a contradiction as before, with one more step needed: since we have ${  g'(t)\over 1+g(g(t))}=1$ and $g(t)\ge t$ for any $t\ge0$, we also have, for any $x\ge0$
$$x=\int_0^{x}{  g'(t)dt\over 1+g(g(t))}=\int_0^{g(x)}{  dt\over 1+g(t)}\le\int_0^{g(x)}{  dt\over 1+t}=\log(1+g(x)) ,$$
whence $e^x\le 1+ g(x)$; if we plug this into the latter inequalities again, we get
$$x=\int_0^{g(x)}{  dt\over 1+g(t)}\le \int_0^{g(x)}e^{-t}dt\le 1 ,$$
as before. By comparison, the same conclusion also holds for $g'=F( {g\circ g})$  with any $F$ analytic and totally monotonic on $(-\epsilon,+\infty)$, and with $F(0)=1$.
\end{remark}
\section{Analysis and discussion:}

\textbf{Numerical solution of $f'=\ Exp ({f^{-1}})$.}

We may take examples for numerical solution of $f'=\exp{f^{ -1}}$ with $f^{-1}(x)=x$ using Rung Kutta method .The disired ODE can be transformed into second order ODE as:
\begin{equation*}\label{ODE}
    v''[t]*v'[t]-v'[t]=0,v(0)=0,v'(0)=2
\end{equation*}

Runge.Kutta methods belong to the class of one step integrators for the numerical solution of ordinary differential equations.In numerical analysis, the Runge.Kutta methods  are a family of implicit and explicit iterative methods, which include the Euler method, used in temporal discretization for the approximate solutions of simultaneous nonlinear equations.\cite{Paul:11} These methods were developed around 1900 by the German mathematicians Carl Runge and Wilhelm Kutta.The most widely known member of the Runge .Kutta family is generally referred to as RK4 , the "classic Runge.Kutta method" or simply as "the Runge.Kutta method''. 

The built in Explicit Runge Kutta method of order 4 is an "embedded" method (i.e., with an embedded error estimation method). Here is the classical method (with the following code )
\begin{verbatim}
    ClassicalRungeKuttaCoefficients[4, prec_] := 
  With[{amat = {{1/2}, {0, 1/2}, {0, 0, 1}}, 
    bvec = {1/6, 1/3, 1/3, 1/6}, cvec = {1/2, 1/2, 1}}, 
   N[{amat, bvec, cvec}, prec]];
\end{verbatim}
Here is a comparison with the built-in method:
\begin{verbatim}
    NDSolve`EmbeddedExplicitRungeKuttaCoefficients[4, MachinePrecision]
ClassicalRungeKuttaCoefficients[4, MachinePrecision]
\end{verbatim}
Now Here is a way to get the numerical solution steps:
\begin{verbatim}
    vf = v /. 
   First@NDSolve[{v''[t]*v'[t] - v'[t] == 0, v[0] == 0, 
      v'[0] == 2}, {v}, {t, 0, 1}, 
     Method -> {"FixedStep", 
       Method -> {"ExplicitRungeKutta", "DifferenceOrder" -> 4, 
         "Coefficients" -> ClassicalRungeKuttaCoefficients}}, 
     StartingStepSize -> 1/10];

Transpose@Flatten[vf[{"Coordinates", {"ValuesOnGrid"}}], 1]
TableForm[%, TableHeadings -> {Range[0, 10], {t, v}}]
\end{verbatim}
The solution are listed in the following table:
$$
\begin{array}{ccc}
  & t & v \\
 0 & 0. & 0. \\
 1 & 0.1 & 0.205 \\
 2 & 0.2 & 0.42 \\
 3 & 0.3 & 0.645 \\
 4 & 0.4 & 0.88 \\
 5 & 0.5 & 1.125 \\
 6 & 0.6 & 1.38 \\
 7 & 0.7 & 1.645 \\
 8 & 0.8 & 1.92 \\
 9 & 0.9 & 2.205 \\
 10 & 1. & 2.5 \\
\end{array}
$$
Here is the plot points (see Figure \ref{fig:my_label} of solution,seems linear:
\begin{figure}[h!]
    \centering
    \includegraphics[width=0.3\textwidth]{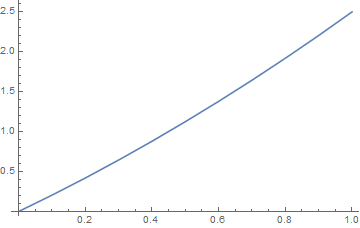}
    \caption{Numerical solutions of     $v''[t]*v'[t]-v'[t]=0,v(0)=0,v'(0)=2$ using RK4 method}
    \label{fig:my_label}
\end{figure}
In the last of this section we may give Pade approximant of the formel power series and ChebychevT (Chebychev function of the first kind) of the obtained formel power series for $n=3$,The obtained pade approximant of the formal power series the solution of (\ref{1}) is rational polynomial with rational coefficients and we have got interesting formula for ChebychevT which is expressed with iteration of $\pi $.
Here is the mathematica code for Pade approximant with $n=3$
\begin{verbatim}
    n = 3;
ansatz = x + Sum[C[k] x^k, {k, 2, n}] + O[x]^(n + 1)
ansatz /. First[SolveAlways[D[ansatz, {x, 2}]
ComposeSeries[D[ansatz, x], InverseSeries[ansatz]] - D[ansatz, x] == 0, x]]
PadeApproximant[Normal[ansatz],{x,0,6]}]

ChebycheVT[R,{x,0,1}]
\end{verbatim}
The obtained Pade Approximant is:
\begin{align*}
   1 +  \frac{253755882643535627x}
  {36875231729969238} + \frac{13180397651429219575x^2}
  {811255098059323236} + \frac{36596037023196361069x^3}{2433765294177969708}&\\
\end{align*}

The obtained ChebychevT in the range $(0,1)$ is given by :
\begin{align*}
    1-\frac{\pi^2 x^2}{8}-\frac{\pi^2 x^3}{8}+\frac{(\pi^4-12\pi^2)x^4}{384}+O[x^5]
\end{align*}

\section{Acknowledgement:}
The author would like to thank Todd Copland  and pietro Majer for their helpful to improve the quality of this paper .

\bibliographystyle{elsarticle harv}








\end{document}